\documentclass[12pt]{amsart}
\usepackage{amsmath,amsfonts,euscript,amscd,amsthm,amssymb,upref,graphics}
\usepackage[retainorgcmds]{IEEEtrantools}

\theoremstyle{plain}
\newtheorem{theorem}{Theorem}
\swapnumbers

\newtheorem{proposition}[subsection]{Proposition}
\newtheorem{lemma}[subsection]{Lemma}

\theoremstyle{definition}
\newtheorem{definition}[subsection]{Definition}

\newtheorem{remark}[subsection]{Remark}

\newtheorem{nothing*}[subsection]{}

\newcommand{\rien}[1]{}

\newcommand{\Spec}{ \operatorname{{\rm Spec}}}

\newcommand{\VFH}{ \operatorname{{\rm VF}_{hol}}}

\newcommand{\LieH}{ \operatorname{{\rm Lie}_{hol}}}
\newcommand{\VFA}{ \operatorname{{\rm VF}_{alg}}}
\newcommand{\LieA}{ \operatorname{{\rm Lie}_{alg}}}

\newcommand{\Aut}{ \operatorname{{\rm Aut}}}

\newcommand{\C}{\ensuremath{\mathbb{C}}}

\newcommand{\nn}{\nonumber}

\def\auta{{\rm Aut}_{alg}}

\renewcommand{\epsilon}{\varepsilon}
\renewcommand{\phi}{\varphi}

\begin{document}
\title[Algebraic density property of Danilov-Gizatullin surfaces  ] 
{Algebraic density property of Danilov-Gizatullin surfaces}
\author{Fabrizio Donzelli}
\address{Institute of Mathematical Sciences, Stony Brook University, 
Stony Brook, NY, 11794 \ \ USA}
\email{fabrizio@math.sunysb.edu}

{\renewcommand{\thefootnote}{} \footnotetext{2000
\textit{Mathematics Subject Classification.} Primary: 32M05,14R20.
Secondary: 14R10, 32M25.}}
\begin{abstract}
A Danilov-Gizatullin surface is an affine surface $V$ which is the complement of an ample section $S$ of a Hirzebruch surface. The remarkable theorem of Danilov and Gizatullin states that the isomorphism class of $V$ depends only on the self-intersection number $S^2$.  
In this paper we apply their theorem to present $V$ as the quotient of an affine threefold by a torus action, and to 
prove that the Lie algebra generated by the complete algebraic vector fields on $V$ coincides with the set of all algebraic vector fields.

\end{abstract}
\maketitle \vfuzz=2pt

\section{Introduction}

\subsection*{Anders\'en-Lempert theory and density property}

\hspace{1cm}

Let $X$ be a complex manifold. We say that a holomorphic vector field $\mu$ on $X$ is complete if the solution of the associated first order ODE exists
for all complex time, for any choice of initial point on $X$.
In other words, $\mu$ induces a holomorphic action of the additive group $\C_+$ on $X$; conversely, given an holomorphic $\C_+$-action $\psi_t$ on $X$
there is a unique complete holomorphic vector field $\mu$ such that its flow coincides with $\psi_t$.

For $n\geq 2$,  the abundance (in the sense of the next definition) 
of complete vector fields on $\C^n$ was a crucial observation 
in the work of Anders\' en and Lempert on the group $\Aut_{hol}(\C^n)$ of holomorphic automorphisms of $\C^n$ \cite{AL}, further developed by Forstneric and Rosay \cite{FR}.

\begin{definition}\label{DP}

A complex manifold $X$ is said to have the density property if 
the Lie algebra $\LieH (X)$ generated by
complete holomorphic vector fields on $X$ is dense, with respect to the compact-open topology, in
the Lie algebra $\VFH (X)$ of all holomorphic vector fields on
$X$. 
\end{definition}

The density property for $\C^n$ has many implications.
An example is the existence of non-equivalent holomorphic embeddings of $\C^k$ into $\C^n$ \cite{BK}, a result that was used by Kutzschebauch and Derksen to construct  non-linearizable holomorphic actions of reductive Lie groups on affine spaces \cite{DK1, DK2}.

The density property was introduced by Varolin \cite{V1},
who was the first to extend the results of Anders\'en and Lempert  
to Stein manifolds different from $\C^n$. An interesting application
from the prospective of the present paper concerns $k$-homogeneity. 

\begin{definition}\label{tr}
Let $X$ be a complex manifold, $k$ a positive integer.
 We say that a group $G$  of holomorphic automorphisms of $X$
acts $k$-transitively on $X$ if given two collections $(x_1,...,x_k)$, $(y_1,...,y_k)$ of pairwise distinct points, there exists an automorphism $f\in G$ such that $f(x_i)=y_i$.
\end{definition}

\begin{proposition}\label{dptr}
If a Stein manifold $X$ has the density property, then the group
$\Aut_{hol}X$ of holomorphic automorphisms of $X$ is acts $k$-transitively for all $k$.
\end{proposition}

Varolin found various examples of Stein manifolds with the density property, 
such as some homogeneous spaces of complex semisimple Lie groups 
(see the paper with Toth, \cite{V2}).

If we are interested in affine algebraic varieties, we can refine the definition as follows.

\begin{definition}\label{ADP}

A smooth affine algebraic variety $X$ is said to have the algebraic density property if 
the Lie algebra $\LieA (X)$ generated by
complete algebraic vector fields on $X$ coincides with 
the Lie algebra $\VFA (X)$ of all  algebraic vector fields on
$X$. 
\end{definition}

Since the ring of regular functions of an affine variety is dense in the ring
of the holomorphic functions, the algebraic density property implies the density property. 
Kaliman and  Kutzschebauch proved that all complex linear algebraic groups with the exception of the tori $(\C^*)^n$ and $\C_+$ have the algebraic density property \cite{KK1}. The tori do not have the algebraic density property since Anders\'en \cite{A}, by using some results from Nevanlinna theory, proves that all complete polynomial vector fields on $(\C^*)^n$ must have divergence zero with respect to the canonical volume form, but we do not know  if the tori have the (not necessarily algebraic) density property.  
Dvorsky, Kaliman and the author \cite{DDK} prove the algebraic density property for homogeneous spaces  not isomorphic to $(\C^*)^n$, of dimension at least three that are quotient of linear algebraic groups by a reductive subgroup. 
The proof relies on the methods developed in \cite{KK1}, but it also requires some non trivial facts from Lie theory,  and an application of the Luna's slice theorem. 
Finally, we mention the affine surface $X$ in $\C^3$ given by the equation $x+y+xyz=1$. By using results of Brunella on holomorphic foliations on rational surfaces \cite{Br1, Br2}, Kaliman and Kutzschebauch \cite{KK3} prove that a complete vector field on $X$ must have divergence zero with respect to the volume form $\frac{dx}{x}\wedge\frac{dy}{y}$, 
and hence that $X$ does not have the density property.
In the same preprint, the authors show that $\Aut_{hol} (X)$
acts $k$-transitively for all $k$. Since we do not know if $X$ has the density property, $X$ could be counterexample to the converse of Proposition \ref{dptr} 
(together with the tori).
For a good survey of Anders\'en-Lempert Theory and it applications, we refer the reader to \cite{KK3}.

\subsection*{Danilov-Gizatullin surfaces}

\hspace{1cm}

A Danilov-Gizatullin surface is an affine surface which is the complement of an ample section $S$ in a Hirzebruch surface $\Sigma_d$. 
The ampleness of $S$ implies that $n>d$ and $n\geq 2$.
We have the following remarkable result of Danilov and Gizatullin \cite{DG}
(see also \cite{FKZ}  for a short proof ):
\begin{theorem}\label{DGT}
Let $S$ ($S'$) be an ample section of $\Sigma_d$ ($\Sigma_{d'}$). 

Then 
 $\Sigma_d-S$ is isomorphic to $\Sigma_{d'}-S'$ if and only if
$S^2=S'^2$. 
\end{theorem}

Consequently, for all $n\geq 2$
we will denote  by $V_n$ the surface complementary to a section $S$ with $S^2=n$.  
The Danilov-Gizatullin surfaces belong to the class of the affine rational
surfaces with trivial Makar-Limanov invariant 
(the suburing of regular functions that are invariant with respect to all 
$\C_+$-actions consists of constants only) \cite{GMMR}. 
Moreover, $V_n$ is a flexible affine variety, that is for any point $x\in V_n$
there are pairs of $\C_+$-actions, such that the corresponding complete vector
fields span the tangent space at $x$. 
On a forthcoming paper of  Arzhantsev, Flenner, Kaliman, Kutzschebauch and Zaidenberg \cite{AFKKZ}, it is proven that the flexibility implies that the subgroup $G$ of $\Aut (V_n)$ generated by $\C_+$-actions acts $k$-transitively on $V_n$, for all $k$ (see Definition \ref{tr} ).
Moreover, the set of complete vector fields whose flow is an algebraic $\C_+$-action generates an infinite-dimensional Lie algebra. 
In this paper (Theorem \ref{main}) we prove that  for all $n\geq 2$, $V_n$ has the algebraic density property.

\subsection*{Sketch of the proof of Theorem \ref{main}}

\hspace{1cm}

By presenting the Hirzebruch 
surface $\Sigma_d$ as a quotient of an open toric variety in $\C^4$ by a two dimensional torus $T^2$, we can describe any section $S$ of the ruling 
as the zero locus of a $T^2$-invariant polynomial in $\C^4$.
 The Danilov-Gizatullin theorem now plays a crucial role:
by choosing to embed $V_n$ in $\Sigma_{n-2}$, and by making a specific choice for a section $S$ with $S^2=n$, we prove that $V_n$ is isomorphic to the algebraic quotient by a one-dimensional torus $T$ of a smooth affine threefold $F_n$
(Theorem \ref{quotient}), and we find a set of generators for
$\C [F_n]^T$, the ring of $T$-invariant regular functions on $F_n$ (Proposition \ref{ring}), which is isomorphic to $\C [V_n]$. Next, we construct some complete $T$-invariant vector fields on $F_n$, that descends to complete vector fieds on $V_n$ (Proposition \ref{fields}). We then perform some computations, involving those fields and the generators of $\C [F_n]^T\cong\C [V_n]$, until we construct a non-zero $\C [V_n]$-submodule $\mathfrak{N}$ of all algebraic vector fields that is contained in the Lie algebra generated by the complete ones (Theorem \ref{M}). We then 
extend $\mathfrak{N}$ to a $\C [V_n]$-submodule $\mathfrak{M}$  such that the fiber of $\mathfrak{M}$ at a point
$x$ of $V_n$ generates the tangent space at $x$.
The transitivity of the action of the group
$\auta{V_n}$ of algebraic automorpshisms of $V_n$ allows then to apply to $\mathfrak{M}$ a technical principle (Theorem \ref{principle}) to prove  the algebraic density property
of $V_n$ (Theorem \ref{main}).

\section{Construction of  $V_n$ as a quotient by a torus action}

A Hirzebruch surface can be realized as a quotient
of a toric variety by a two-dimensional torus (\cite{K}, Section 1):
consider the action of  $T^2\cong\C_{t_1}^*\times \C_{t_2}^*$ on $\C^4_{(a_1,a_2,a_3,a_4)}$ given by 

\begin{align}
\nonumber (t_1,t_2).(a_1,a_2,a_3,a_4)=(t_1t_2^da_1,t_2a_2,t_1a_3,t_2a_4)
\end{align}

and let

\begin{align}
\nn Z=\{a_1=a_3=0\}\cup \{a_2=a_4=0\}.
\end{align}

Then we have that 

\begin{align}
\nn \Sigma_d\cong (\C^4-Z)/T^2. 
\end{align}

The isomorphism above can be understood by defining coordinates

\begin{align}
&\nn t_0=\frac{a_1}{a_2^da_3}\quad \text{for}\quad a_2\neq 0 ,\\
&\nn t_\infty=\frac{a_1}{a_4^da_3}\quad\text{for}\quad a_4\neq 0 ,\\
&\nn v=\frac{a_2}{a_4}\in\mathbb{P}^1
\end{align}

(since they are $T^2$-invariant, they are  well defined).

The quotient space then is isomorphic to the Hirzerbruch surface $\Sigma_d$, 
with ruling
 
\begin{align}
&\nn \pi : \Sigma_d\rightarrow \mathbb{P}^1\\ 
&\nonumber \pi (a_1,a_2,a_3,a_4)=(a_2,a_4)
\end{align}

and transition function

\begin{align}
\nn t_\infty=v^dt_0.
\end{align}

Recall now (\cite{B}, chapter IV)
 that the Picard group of $\Sigma_d$ is generated by the linear equivalence class $F$ of the fiber of $\pi$, and by the class of the unique irreducible curve $C$ with self-intersection $C^2=-d$.
It follows that  $S$ is linear equivalent to $C+bF$, and that $n=2b-d$.

Theorem \ref{DGT} implies that $\Sigma_d-S\cong V_n$ for any section $S$ with $S^2=n$, and we claim that

\begin{align}
\nn t_\infty=v^b 
\end{align}

is a local equation (for $v\neq \infty$) of a section
$S$ of the ruling $\pi$ with self intersection $n=2b-d$. In fact, the base locus of the
family of section $t_\infty=cv^b$ consists of two multiple points, namely
the point $A(v=0,t_\infty=0)$ with multiplicity $b$, and $B(v=\infty, t_0=0)$ with multiplicity $b-d$: if we sum the intersection multiplicities at $A$ and 
$B$, we obtain that $S^2=2b-d=n$.

In the coordinates $(a_1,a_2,a_3,a_4)$ the defining equation of $S$ takes the
form

\begin{align}\label{a1}
a_1a_4^{b-d}-a_2^ba_3=0
\end{align}
(the ampleness of $S$ implies that $b>d$).

Let $S'$ be the hypersurface in $\C^4$ given by equation (\ref{a1}).
Then the inverse image $\rho^{-1}(S)=S'-Z$ under the quotient map is the closed 
$T$-invariant subset of $\C^4$ given by the equation (\ref{a1}),  and

\begin{align}\label{a2}
\nn V_n\cong\rho (\C^4-(S'\cup Z))=(\C^4-(S'\cup Z))/T^2. 
\end{align}

Since $Z\subset S'$, the variety $\C^4-(S'\cup Z)$ is isomorphic to the affine manifold
$V'=\C-S'$.
\begin{lemma}\label{gq=aq}
The geometric quotient $V'/T^2$ is isomorphic to the algebraic quotient

$V'//T^2=\Spec \C [V']^{T^2}$, where $\C [V']^{T^2}$ denotes the subring of the $T^2$-invariant regular functions of $V'$.

\begin{proof}
The orbits of the $T^2$-action restricted to $V'$ are all closed. 
Therefore the lemma follows  from the Luna's slice theorem (\cite{Dr}, Theorem 5.4).
\end{proof}

\end{lemma}

We can invoke Theorem \ref{DGT} again and fix 
from now on that $V_n$ is embedded in $\Sigma_{n-2}$, that is $n=d+2$, and $b=d+1$. Under this assumption, the equation of $S'$ becomes

\begin{align}
a_1a_4-a_2^ba_3=0.
\end{align}

\begin{proposition}
Consider the smooth affine threefold 

\begin{align}\label{equation}
F_n=\{ a_1a_4-a_2^ba_3=1 \}\subset \C^4.
\end{align}

Then $V'$ is equivarianlty isomorphic to 
$\C^*\times F_n$, where the action of $T^2$ on $\C^*\times F_n$ is defined by:

\begin{align}\label{action}
(t_1,t_2).(w,a_1,a_2,a_3,a_4)=
(t_1^{-1}t_2^{-b}w,t_2^{-1}a_1,t_2a_2,t_2^{-b}a_3,t_2a_4).
\end{align}

\begin{proof}
The isomorphism (as affine varieties)
is given by

\begin{align}
(a_1,a_2,a_3,a_4)\mapsto \left( a_1a_4-a_3^ba_2,\frac{a_1}{a_1a_4-a_3^ba_2}, a_2,\frac{a_3}{a_1a_4-a_3^ba_2},a_4 \right)
\end{align}

Then it is easy to check that the induced $T^2$ action on $\C^*\times F_n$ by the isomorphism is given by (\ref{action}).

\end{proof}

\end{proposition}

When we pass to the quotient, we get rid of the factor $\C^*$, as follows.

\begin{theorem}\label{quotient}
For any $n\geq 2$, the Danilov-Gizatullin surface $V_n$ is isomorphic to the algebraic quotient $F_n//T$, where $F_n$ is the affine threefold given
by the equation 
\begin{align}
a_1a_4-a_2^ba_3=1
\end{align}
 and the torus $T\cong \C^*$ acts on $F_n$ via

\begin{align}
t.(a_1, a_2,a_3,a_4)= (t^{-1}a_1, ta_2, t^{-b}a_3,ta_4)
\end{align}

\begin{proof}
Lemma \ref{gq=aq} implies that $V'/T^2\cong V'//T^2$.
Observe that the $T^2$-action of $V'$ (equation \ref{action}) is such that 
the $t_1$-variable acts non-trivially only on the $w$-coordinate.
Therefore there are no non-constant $T^2$-invariant polynomials depending
on $w$, and the result follows.
\end{proof}

\end{theorem}

Perhaps, the above theorem has been already estabilished, but we do not know a reference for it.

\begin{remark}
Theorem \ref{quotient} presents $V_n$ as a quotient of a hypersurface  in 
$\C^4$ given by an equation of the type $uv=p(y,z)$, for $p$ a polynomial
with smooth zero locus. The affine manifolds of this type have  
the algebraic density property, as shown by Kaliman and Kutzschebauch \cite{KK2}. Furthermore, it is clear from its defining equation that 
$F$ is a ramified $b$-sheeteed covering of $SL_2$, another manifold with the algebraic density property \cite{V2}. 
The behavior of the density property with respect to 
quotients or coverings is not known in general.
\end{remark}

\begin{proposition}\label{ring}
$\C [V_n]\cong \C [F_n]^T$ is generated by the monomials 

\begin{align}
y= a_1a_2,\quad z=a_1a_4,\quad \quad
 x_k=a_2^{b-k}a_3a_4^k \quad (0\leq k\leq b). 
\end{align}

\begin{proof}
Since the action is monomial, it is sufficient to look
for monomial generators. The action of $t$ on $a_1^Xa_2^Ya_3^Za_4^W$ is 

\begin{align}
\nn t.(a_1^Xa_2^Ya_3^Za_4^W)=t^{-X+Y-bZ+W}a_1^Xa_2^Ya_3^Za_4^W,
\end{align}

therefore a monomial is invariant if and only if 

\begin{align}\label{lattice}
-X+Y-bZ+W=0.
\end{align}

We show that every solution (with non-negative integers) of equation (\ref{lattice})
is a linear combination with non-negative integer coefficients of the vectors

\begin{align}
&\nn (1,1,0,0)\quad (1,0,0,1),\\
&\nn (0,b-k,1,k)\quad (\text{for}\quad 0\leq k\leq b).
\end{align}

Write $W=NbZ+r$ and $Y=MbZ+r'$, for  $0\leq r,r'\leq bZ$.
Suppose that $W<bZ$: choose $Z$ integers $k_i$, with $0\leq k_i\leq b$,
such that $\sum_i^Zk_i=W$. Then

\begin{align}
\nn (X,Y,Z,W)=X(1,1,0,0)+ \sum_i^Z(0,b-k_i,1,k_i)
\end{align}

Suppose instead that $W\geq bZ$. Then $X-MbZ=W-bZ+r'\geq 0$; choose $Z$ integers $k_i$, with $0\leq k_i\leq b$, such that $\sum _i^Z(b-k_i)=r'$. Then we can write

\begin{align}
\nn (X,Y,Z,W)=(X-MbZ)(1,0,0,1)+\sum_i^Z(0,b-k_i,1,k_i)+MbZ(1,1,0,0).
\end{align}

\end{proof}

\end{proposition}

\section{Algebraic density property of $V_n$}

In this section we prove the main result of the paper. 

\begin{theorem}\label{main}
The Danilov-Gizatullin surfaces $V_n$ have the algebraic density property.
\end{theorem}

We make use of the following principle \cite{KK1}. 

\begin{theorem}\label{principle}
Let $X$ be a smooth affine algebraic variety, such that the group $\auta (X)$
of algebraic automorphisms acts transitively on it. Let $\mathfrak{M}$ be a submodule
of the $\C [X]$-module of all algebraic vector fields such that $\mathfrak{M}\subset\LieA (X)$. Suppose that the fiber of $\mathfrak{M}$ 
at a point $x_0\in X$ generates the tangent space $T_{x_0}M$.
Then $X$ has the algebraic density property.
\end{theorem} 
 
In our case the transitivity of the action of $\auta (V_n)$ follows from the following result of Gizatullin \cite{G}.

\begin{theorem}
Let $X$ be an irreducible smooth affine variety over an algebraically closed field of characteristic zero. If $X$ can be completed by a smooth rational curve, then $\auta (X)$ acts transitively on $X$.
\end{theorem}

We start by listing the complete vector fields that will be used to construct the module $\mathfrak{M}$ of the theorem above.

\begin{proposition}\label{fields}
The following vector fields on $F_n$ are complete and $T$-invariant, and they descend to non-identically zero complete vector fields on $V_n$.

\begin{align}
&\delta =ba_2^{b-1}a_3\frac{\partial}{\partial a_1}+a_4\frac{\partial}{\partial a_2}\\
&\delta'=a_1^{b-1}a_2^b\frac{\partial}{\partial a_4}+a_1^b\frac{\partial}{\partial a_3}\\
&\epsilon=a_1\frac{\partial}{\partial a_1}-a_4\frac{\partial}{\partial a_4}
\end{align}

Moreover, $\delta$ and $\delta'$ are locally nilpotent derivations
(their flow is an algebraic action of $\C_+$).

\begin{proof}
The vector fields are tangent to $F_n$, since they annihilate the defining equation of the threefold. The invariance and their completeness is a straightforward check left to the reader. 
Finally, the vector fields descends non-trivially on $V_n$ because they are generically transversal to the vector field $
-a_1\frac{\partial}{\partial a_1}+a_2\frac{\partial}{\partial a_2}
-ba_3\frac{\partial}{\partial a_3}+a_4\frac{\partial}{\partial a_4}
$ that defines the $T$-action on $F_n$.

\end{proof}

\end{proposition}

We collect in the next two lemmas the formulas upon which the 
proofs of lemmas 4.5-4.9 are based.

\begin{lemma}\label{on functions}
For the surface $V_n$:
\begin{align}
&\epsilon (x_k)=-kx_k\\
&\epsilon (y)=y\\
&\delta (y)=1+nx_0\\
&\delta (x_k)=(b-k)x_{k+1}\\
&\delta' (x_0)=y^b\\
&\delta'(y)=0
\end{align}
\end{lemma}

\begin{lemma}\label{commutation}
For the surface $V_n$:

\begin{align}
&[\epsilon ,\delta ] = -\delta\\
&[\epsilon , \delta '] = b\delta' \quad (n=b+1) 
\end{align}

\end{lemma}

We will also need a well-known fact (\cite{V4}, Proposition 3.1) about differential equations that will be used implicitly in some of the remaining statements.

\begin{proposition}
Let $\mu$ be a complete vector field on a complex manifold $X$, and $f$ be a holomorphic function.
Then the vector field $f\mu$ is complete if and only if $\mu^2(f)=0$.
\end{proposition}

\begin{lemma}\label{epsilon}
The following algebraic vector fields belong to $\LieA (V_n)$:

(1) $x_0\epsilon$, (2) $x_b\epsilon$, (3) $x_0x_b\epsilon$.

\begin{proof}
Statement (1) follows from $\epsilon (x_0)=0$ (in particular, $x_0\epsilon$ is complete). 

As for (2), let $X_1= [\delta , x_0\epsilon]$, $X_s= [\delta , X_{s-1}]$.
By definition $X_s\in\LieA (V_n)$, for all $s$.
Then by induction on $s\geq 2$ it follows that
\begin{align}
\nn X_s=sb(b-1)...(b-s+2)x_{s-1}\delta +b(b-1)...(b-s+1)x_s\epsilon.
\end{align}

In particular, for $s=b$, one obtains 
$X_b=b!x_b\epsilon + b^2(b-1)...(2)x_{b-1}\delta$.
Since $\delta^2(x_{b-1})=0$, $x_{b-1}\delta$ is complete, and the statement (2)
is proven. 

From the first two facts, one has that (3) $[x_0\epsilon, x_b\epsilon ]=-bx_0x_b\epsilon\in\LieA (V_n)$. 

\end{proof}

\end{lemma}

\begin{lemma}\label{delta}
For $0\leq k \leq b$, for all $N>0$, $x_kx_b^N\delta\in\LieA (V_n)$.

\begin{proof}

Now  we calculate that 

\begin{align}
\label{s1}
[x_0\epsilon , x_b\delta ]&=-(1+b)x_0x_b\delta-bx_1x_b\epsilon\\
\label{s2}[\delta , x_0x_b\epsilon ]&=x_0x_b\delta+bx_1x_b\epsilon
\end{align} 

The sum of  (\ref{s1}) and (\ref{s2}) shows  that
$x_0x_b\delta\in\LieA (V_n)$ (use Lemma \ref{epsilon}).
Then, for all $M\geq 0,$ we have $[x_b^M\delta , x_0x_b\delta ]=bx_1x_b^{M+1}\delta\in\LieA (V_n)$.
By induction on $k$ we then prove the statement for all positive $k$.

As for $k=0$, observe that since $yx_b=(1+x_0)x_{b-1}$,
we obtain $yx_b\delta =x_0x_{b-1}\delta + x_{b-1}\delta$.
This equation shows that $yx_b\delta\in\LieA (V_n)$: in fact, $x_{b-1}\delta$
is complete, and $x_0x_{b-1}\delta\in\LieA (V_n)$ as follows from Lemma \ref{epsilon} and Lemma \ref{delta} applied to

\begin{align}
\nn [x_0\epsilon , x_{b-1}\delta ]=-bx_0x_{b-1}\delta-bx_1x_{b-1}\epsilon = -bx_0x_{b-1}\delta-bx_0x_b\epsilon, 
\end{align}
the second equality being true since $x_1x_{b-1}=x_0x_b$.
Hence, for $M\geq 0$ we obtain

\begin{align}
\nn [x_b^M\delta , yx_b\delta ]=x_b^{M+1}(1+nx_0)\delta =x_b^{M+1}\delta+nx_0x_b^{M+1}\delta
\end{align} 

which shows the case $k=0$.

\end{proof}

\end{lemma}

\begin{lemma}\label{xe}
For $1\leq k \leq b$, $N>0$, $M\geq 0$,  $x_0^Mx_kx_b^N\epsilon\in\LieA (V_n)$.

\begin{proof}
We know that $x_0x_b\epsilon\in\LieA (V_n)$ (Lemma \ref{epsilon}).

To show that  $x_kx_b\epsilon\in\LieA (V_n)$ for all $k$, apply induction on $k$
and Lemma \ref{delta} to the equation

\begin{align}
\nn [\delta ,x_{k-1}x_b\epsilon ]=x_{k-1}x_b\delta + (b-k+1)x_kx_b\epsilon.
\end{align}

Next, the result for $N'\geq 0$ $M=0$ follows from the equation
\begin{align} 
\nn [x_b^{N'}\delta , x_kx_b\epsilon ]= (1+bN')x_kx_b^{N'+1}\delta+x_b^{N'+1}(b-k)x_{k+1}\epsilon.
\end{align}

Finally, the case $M\geq 0$, $N>0$ follows from 

\begin{align}
\nn [x_0^M\epsilon , x_kx_b^N\epsilon ]=(-Nb-k)x_0^Mx_kx_b^N\epsilon.
 \end{align}

\end{proof}
\end{lemma}

The next computations deal with the $y$-variable.

\begin{lemma}\label{ey}
For all $R\geq 0$,
$y^{b+R}\epsilon\in\LieA (V_n)$.
\begin{proof}

Use the formula

\begin{align}
\nn [x_0\epsilon , y^R\delta' ]= (b+j)x_0y^R\delta'-y^{b+R}\epsilon
\end{align}
and notice that $\delta'^{2}(y^Rx_0)=0$.
\end{proof}

\end{lemma}

\begin{proposition}\label{d}
For $N>0$, $M>0$, $0<k<b$, $R\geq 0$,
$y^{b+R}x_0^Mx_kx_b^N\epsilon\in\LieA (V_n)$

\begin{proof}

Apply Lemma \ref{xe} and Lemma \ref{ey} to the equation

\begin{align}
&\nn [y^{b+R}\epsilon ,x_0^Mx_kx_b^N\epsilon ]=\\
&\nn (1-k-bN-b-R) y^{b+R}x_0^Mx_kx_b^N\epsilon.
\end{align}

\end{proof}

\end{proposition}

Next, we make use of some relations between the invariant monomials.

\begin{proposition}\label{monomial}
For any collection $k_1,...,k_r$, there are integers $N>0$ and $M>0$,
$0\leq h\leq b$, such that

\begin{align}
x_0x_{k_1}x_{k_2}...x_{k_r}x_b=x_0^Mx_hx_b^N,
\end{align}

\begin{proof}
We observe that:
\begin{align}
&if \quad h+k\leq b\quad x_kx_h=x_{h+k}x_0\label{red1}\\
&if \quad h+k\geq b\quad x_kx_h=x_bx_{h+k-b}\label{red2} 
\end{align}
Write the monomial in the form

\begin{align}
\nn x_0^{N'}x_{k_1}...x_{k_s}x_b^{M'}
\end{align}
for $x_{k_1},...,x_{k_{s}}\neq x_0,x_b$.

Then, by applying (\ref{red1}) and (\ref{red2})
we can eliminate the factors $x_k$ not equal to $x_0$ and $x_b$
untill there will be left at most one.

\end{proof}

\end{proposition}

\begin{theorem}\label{M}
Let $\mathfrak{N}$ be the
$\C [V_n]$-submodule of $\VFA (V_n)$ generated by the vector field
$x_0x_1...x_by^b\epsilon$. Then $\mathfrak{N}\subset\LieA (V_n)$.
\begin{proof}
We need to show that, for all monomials $z^Ty^Rx_{k_1}...x_{k_r}$,

$z^Ty^Rx_{k_1}...x_{k_r}x_0x_1...x_by^b\epsilon\in\LieA (V_n)$.
Observe first that we can assume $T=0$, since we can eliminate $z$ via the relation $z=1+x_0$.

Then, according to Proposition \ref{monomial}, we can write

\begin{align}
\nn y^Rx_{k_1}...x_{k_r}x_0x_1...x_by^b=y^{b+R}x_0x_{h_1}...x_{h_{r'}}x_b=
\nn x_0^Nx_hx_b^M,
\end{align}
for some $N,M>0$, $0\leq h\leq b$

The theorem follows from Proposition \ref{d}.

\end{proof}

\end{theorem}

\subsection*{Conclusion of the proof of Theorem \ref{main}}
We follow the idea of \cite{KK1} to produce a module $\mathfrak{M}$
satisfying the hypothesis of Theorem \ref{principle}.
Consider the regular function $f=x_b$.
The flow of the vector field $f\delta$ is an algebraic action of $\C_+$:
this can be checked directly by solving the associated ODE, or
by using the fact that $\delta$ is a locally nilpotent derivation and 
$\delta (x_b)=0$ (\cite {Fre}, Principle 7 page 24).
Let $p\in V_n$ be  a point with $f(p)=0$, $\delta_p\neq 0$, 
$\epsilon_p\neq 0$, $\epsilon_p (f)\neq 0$:
then (\cite{KK1}, Claim page 4) the flow $\phi$ of $f\delta$ at time one
 induces an isomorphism $\phi_*$
on the tangent space at $p$, 
that map $\epsilon_p$ to $\epsilon_p+\epsilon_p(f)\delta_p$.
Let $\mathfrak{M}=\mathfrak{N}+\phi_*\mathfrak{N}$: then $\mathfrak{M}\subset\LieA (V_n)$, and since $\epsilon_p$ and $\delta_p$ are linearly independent vectors,
its fiber at $p$ spans the tangent space at $p$. 
Since $V_n$ is homogeneous with respect to $\auta (V_n)$, 
we can then apply Theorem \ref{principle} to conclude the proof of Theorem \ref{main}.

\subsection*{Acknowledgments}
I would like to thank Shulim Kaliman and Dror Varolin 
for reading the paper and giving a lot of useful comments; 
Patrick Clarke and Andrew Young for useful conversations. 

\vfuzz=2pt

\providecommand{\bysame}{\leavevmode\hboxto3em{\hrulefill}\thinspace}

\end{document}